\documentclass[12pt, a4paper]{article}
\usepackage{float}
\usepackage{titlesec} %自定义多级标题格式的宏包
\usepackage{geometry}
\usepackage{multirow}
\usepackage{enumerate}
\usepackage{booktabs}
\usepackage{ytableau}
\usepackage{tikz}
\usepackage{graphicx}
\usepackage{caption}
\usepackage{tikz}
\usepackage{mathrsfs}
\usepackage{makecell}
\usepackage{algorithm,algorithmic}
\usepackage{lmodern}
\usepackage{amssymb}
\usepackage{amsfonts}
\usepackage{graphicx}
\usepackage{epstopdf}
\usepackage{setspace}
\usepackage[unicode={ture},colorlinks,
            linkcolor=black,
            anchorcolor=blue,
            citecolor=green]{hyperref}
\usepackage{appendix}
\usepackage{amsmath}
\usepackage{amsthm}

\usepackage{authblk}

\usepackage{xcolor}
\usepackage[normalem]{ulem}

\allowdisplaybreaks[4]

\newtheorem{thm}{Theorem}[section]
\newtheorem{lem}[thm]{Lemma}
\newtheorem{pr}[thm]{Proposition}
\newtheorem{cor}[thm]{Corollary}

\theoremstyle{definition}
\newtheorem{def1}{Definition}[section]
\newtheorem{re}{Remark}[section]
\newtheorem{prob}{Problem}[section]

\title{Aldous' spectral gap property for normal Cayley graphs on symmetric groups}

\author{Yuxuan Li, Binzhou Xia and Sanming Zhou\\ 
School of Mathematics and Statistics, The University of Melbourne, Parkville, VIC 3010, Australia}
\begin{document}
\maketitle
\openup 0.6 \jot % Added by Sanming Zhou

\renewcommand{\thefootnote}{\empty}%{footnote}}
\footnotetext{E-mail addresses: yuxuan11@student.unimelb.edu.au (Yuxuan Li), binzhoux@unimelb.edu.au (Binzhou Xia), sanming@unimelb.edu.au (Sanming Zhou)}

\begin{abstract}
Aldous' spectral gap conjecture states that the second largest eigenvalue of any connected Cayley graph on the symmetric group $S_n$ with respect to a set of transpositions is achieved by the standard representation of $S_n$. This celebrated conjecture, which was proved in its general form in 2010, has inspired much interest in searching for other families of Cayley graphs on $S_n$ with the property that the largest eigenvalue strictly smaller than the degree is attained by the standard representation of $S_n$. In this paper, we prove three results on normal Cayley graphs on $S_n$ possessing this property for sufficiently large $n$, one of which can be viewed as a generalization of the ``normal'' case of Aldous' spectral gap conjecture.
\end{abstract}

\section{Introduction}
\label{sec:intro}

We only consider finite simple undirected graphs in this paper. Suppose that $\Gamma$ is such a graph and $A(\Gamma)$ is its adjacency~matrix. Since $A(\Gamma)$ is real and symmetric, all its eigenvalues are real numbers, and they are called the \emph{eigenvalues} of $\Gamma$. We always arrange the eigenvalues of $\Gamma$ in non-ascending order as $\lambda_1\geq \lambda_2\geq \cdots\geq \lambda_n$. Whenever we want to stress the dependence of the $i$-th largest eigenvalue of $\Gamma$ or a real symmetric matrix $M$, we write $\lambda_i (\Gamma)$ or $\lambda_i(M)$ in place of $\lambda_i$. It is known that the largest eigenvalue $\lambda_1 (\Gamma)$ of any regular graph $\Gamma$ is equal to the degree of $\Gamma$.

 %The {\bf degree} of vertex $v_i$ is the number of vertices in $\Gamma$ which are adjacent to $v_i$. We say that $\Gamma$ is {\bf $k$-regular} if every vertex in $V(\Gamma)$ has the same degree $k$. Cayley graphs are an important class of regular graphs in both theory and application fields. 

Let $G$ be a finite group with identity element $e$, and let $S$ be an inverse-closed subset of $G\setminus \{e\}$. The \emph{Cayley~graph} on $G$ with respect to $S$, denoted by $\mathrm{Cay}(G, S)$, is the $|S|$-regular graph with vertex set $G$ and edge set $\{\{g,gs\}~|~g\in G, s\in S\}$. It is readily seen that $\mathrm{Cay}(G,S)$ is connected if and only if its \emph{connection set} $S$ is a generating subset of $G$. We say that $\mathrm{Cay}(G,S)$ is a \emph{normal} Cayley graph if $S$ is closed under conjugation.%\footnote{This definition is different from the following notion defined in \cite{WWX}: A Cayley graph $\Gamma=\text{Cay}(G,S)$ is called a normal Cayley graph if the right regular representation of $G$ is normal in $\text{Aut}(\Gamma)$.}.

It is widely known that the representation theory of finite groups (see \cite{I,JL,S1,S}) plays a critical role in determining eigenvalues of Cayley graphs. In what follows we use 
$$
\widehat{G}=\{\rho_1,\rho_2,\ldots,\rho_k\}
$$ 
to denote a complete set of inequivalent (complex) irreducible matrix representations of $G$, with the assumption that $\rho_1$ is the trivial representation. For any $\rho_i \in \widehat{G}$, the map 
$$
\chi_i: g \mapsto \mathrm{Trace}(\rho_i(g)),\;\, g \in G
$$ 
is the \emph{character} of $\rho_i$, and the ratio 
$$
\tilde{\chi}_i(g):=\frac{\chi_i(g)}{\chi_i(e)}
$$ 
is known as the \emph{normalized character} of $\rho_i$ on $g\in G$, where $\chi_i(e)$ equals the dimension $\dim\rho_i$ of $\rho_i$. Note that $\dim \rho_1 = 1$ for the trivial representation $\rho_1$.

It is known \cite{MS} that the adjacency matrix of $\mathrm{Cay}(G, S)$ equals $\sum_{s\in S} R_{\mathrm{reg}}(s)$, where $\mathrm{reg}$ is the right regular representation of $G$ and $R_{\mathrm{reg}}(s)$ is the permutation matrix depicting the multiplication on $G$ from the right by $s$. Set
$$
\rho_i(S) := \sum_{s \in S} \rho_i(s)
$$
and denote by $\oplus$ the direct sum of matrices. It is known that the adjacency matrix of $\mathrm{Cay}(G,S)$ is similar to $d_1 \rho_1(S)\oplus d_2\rho_2(S)\oplus \cdots\oplus d_k\rho_k(S)$ (see \cite[Proposition 7.1]{MS}), where $d_i$ is the dimension of $\rho_i\in \widehat{G}$ and $d_i\rho_i(S)$ is the direct sum of $d_i$ copies of $\rho_i(S)$.
This implies that the multiset of eigenvalues of $\mathrm{Cay}(G,S)$ is the union of $d_i$ multisets of eigenvalues of $\rho_i(S)$, for $1\leq i\leq k$. Since $\mathrm{Cay}(G,S)$ is regular with degree $|S|$, its largest eigenvalue is equal to $|S| = d_1\rho_1(S)$. Moreover, if $\mathrm{Cay}(G,S)$ is connected, then the multiplicity of $|S|$ is $1$. In the case when $\mathrm{Cay}(G,S)$ is normal, by Schur's Lemma, all $\rho_i(S)$'s are scalar matrices (see \cite[Lemma 5]{PS}) and the eigenvalues of $\mathrm{Cay}(G,S)$ can be expressed in terms of the irreducible characters of $G$ in the following way. 

\begin{pr}[\cite{PS,Z}]\label{prop1.2}
Let $\{\chi_1,\chi_2,\ldots,\chi_k\}$ be a complete set of inequivalent irreducible characters of $G$. Then the eigenvalues of any normal Cayley graph $\mathrm{Cay}(G, S)$ on $G$ are given by 
\begin{equation*}
	\lambda_j=\frac{1}{\chi_j(e)}\sum_{s\in S}\chi_j(s)=\sum_{s\in S} \tilde{\chi}_j(s),\quad j=1,2,\ldots,k.
\end{equation*}
Moreover, the multiplicity of $\lambda_j$ is equal to $\sum_{1\leq i\leq k,~\lambda_i=\lambda_j}\chi_i(e)^2.$
\end{pr}

We say that the second largest eigenvalue of a connected Cayley graph $\mathrm{Cay}(G,S)$ is \emph{attained} or \emph{achieved} by the representation $\rho_i$ of $G$ if
$$
\lambda_2(\mathrm{Cay}(G,S))=\lambda_1(\rho_i(S)).
$$
In view of Proposition \ref{prop1.2}, if $\mathrm{Cay}(G,S)$ is connected and normal, then its second largest eigenvalue is achieved by $\rho_i$ if and only if $\lambda_2(\mathrm{Cay}(G,S))=\sum_{s\in S} \tilde{\chi}_i(s)$.
 
The second largest eigenvalue of graphs has attracted much attention over the past more than three decades (see, for instance, \cite{A,N}). In particular, the second largest eigenvalue of Cayley graphs has been a focus of study for a long time owing to the fact that some important expanders are Cayley graphs. See \cite{HLW} for a survey on expander graphs with applications and \cite[Section 8]{LZ} for a collection of results on the second largest eigenvalue of Cayley graphs. One of the most important results about the second largest eigenvalue of Cayley graphs is Aldous' spectral gap conjecture, which was made by Aldous \cite{A1} in 1992 and completely proved by Caputo, Liggett  and Richthammer \cite{CLR} in 2010. As usual, let $S_n$ and $A_n$ be, respectively, the symmetric and alternating groups on $[n]=\{1,2,\ldots,n\}$, where we assume $n\ge 3$. %A \emph{transposition} in $S_n$ is a permutation exchanging exactly two points of $[n]$. 
Recall from the representation theory of symmetric groups \cite{Sagan} that for each partition of $n$ (which will be defined in Section \ref{sec:pre}), we can construct an irreducible representation of $S_n$ known as its Specht module. It is well known that all Specht modules form a complete list $\widehat{S_n}$ of inequivalent irreducible representations of $S_n$. In particular, the Specht module corresponding to the partition $(n-1,1)$ of $n$ is called the \emph{standard representation} of $S_n$ and is denoted by $\rho_{(n-1,1)}$. Aldous \cite{A1} conjectured that for any generating subset $T$ of $S_n$ which consists of transpositions, the second largest eigenvalue of $\mathrm{Cay}(S_n,T)$ is achieved by the standard representation of $S_n$; that is, $\lambda_2(\mathrm{Cay}(S_n,T))=\lambda_1(\rho_{(n-1,1)}(T)).$ Following a number of attempts for special cases (see, for example,  \cite{C,PS,FOW,HJ}), Aldous' conjecture in its general form was finally proved in \cite{CLR} with the help of a nonlinear mapping and a complicated estimation called the Octopus Inequality. (These two key ingredients appeared almost simultaneously in \cite{D1}, where the Octopus Inequality was proved in some special cases.) In 2016, Cesi \cite{C2} presented a simpler and more transparent proof of the Octopus Inequality which makes it possible to look at Aldous' spectral gap conjecture from an algebraic perspective.  
 	
%If $S$ is a generating subset of the alternating group $A_n$, then $\mathrm{Cay}(A_n,S)$ is connected and $\mathrm{Cay}(S_n,S)$ is the union of two copies of $\mathrm{Cay}(A_n,S)$. Thus $\text{Cay}(S_n, S)$ has the same eigenvalues as $\text{Cay}(A_n,S)$, except that the multiplicity of each eigenvalue of $\text{Cay}(S_n, S)$ is twice the multiplicity of the corresponding eigenvalue of $\text{Cay}(A_n,S)$. If \begin{equation*}
%	\lambda_2(\mathrm{Cay}(A_n,S))=\lambda_3(\mathrm{Cay}(S_n,S))=\lambda_1(\rho_{(n-1,1)}(S)),
%\end{equation*}
%then we also say that $\mathrm{Cay}(A_n,S)$ possesses the Aldous property.

A problem closely related to Aldous' spectral gap conjecture is to determine the exact value of the second largest eigenvalues of some connected Cayley graphs on $S_n$ or $A_n$ with connection set not necessarily formed by transpositions only. This problem is quite challenging in general but has been settled for several families of Cayley graphs on symmetric or alternating groups, including the pancake graphs \cite{C1}, the reversal graphs \cite{CT}, a family of graphs which contains all pancake graphs \cite{CT}, three families of Cayley graphs on $A_n$ \cite{HH}, $\mathrm{Cay}(S_n, C(n,k))$ with $k \in \{n-1, n\}$ even \cite{SZ}, and $\mathrm{Cay}(A_n, C(n,k))$ with $k \in \{n-1, n\}$ odd \cite{SZ}, where $n > 4$ and $C(n,k)$ is the set of $k$-cycles in $S_n$. Let $C(n,k; r)$ be the set of $k$-cycles in $S_n$ which move all points in $\{1, 2, \ldots, r\}$, where $1 \le r < k < n$. In \cite{SZ}, Siemons and Zalesski obtained a lower bound on the second largest eigenvalue of $\mathrm{Cay}(G, C(n,k; r))$ and conjectured that this bound is tight, where $G = S_n$ when $k$ is even and $G = A_n$ when $k$ is odd. In the same paper they proved their conjecture for $k = r+1$ and any $1 \le r < n-1$, and earlier results confirmed this conjecture for $(k, r) = (3, 1), (3, 2)$ \cite{HH} and $(k,r)=(2,1)$ \cite{FOW}. Note that $\mathrm{Cay}(G, C(n,k; r))$ is not normal, but the normalizer of $C(n,k; r)$ in $S_n$ is large, namely it is isomorphic to $S_r \times S_{n-r}$.  

Recently, several researchers have started to work on generalizing Aldous' conjecture to Cayley graphs on symmetric groups whose connection sets contain non-transpositions. Before discussing their results, let us first introduce the following definition, in which the strictly second largest eigenvalue is considered as the Cayley graph involved is not required to be connected. In general, the \emph{strictly second largest eigenvalue} of a regular graph is the largest eigenvalue strictly smaller than the degree of the graph.
 
\begin{def1}
\label{def:aldous property}
We say that a Cayley graph $\mathrm{Cay}(S_n,S)$ on $S_n$ has the \emph{Aldous property} if its strictly second largest eigenvalue is attained by the standard representation of $S_n$, that is, 
$$
\lambda_{t+1}(\mathrm{Cay}(S_n,S))=\lambda_1(\rho_{(n-1,1)}(S)),
$$
where $t := [S_n:\langle S \rangle]$ is the index of $\langle S\rangle$ (the subgroup of $S_n$ generated by $S$) in $S_n$.
\end{def1}

Note that the largest eigenvalue $|S|$ of $\mathrm{Cay}(S_n,S)$ has multiplicity $t$ as $\mathrm{Cay}(S_n,S)$ is the union of $t$ copies of the connected Cayley graph $\mathrm{Cay}(\langle S\rangle,S)$ with degree $|S|$. So $\lambda_{t+1}(\mathrm{Cay}(S_n,S))$ is indeed the strictly second largest eigenvalue of $\mathrm{Cay}(S_n,S)$. In view of Proposition \ref{prop1.2}, when $\mathrm{Cay}(S_n,S)$ is normal, it has the Aldous property if and only if
$$
\lambda_{t+1}(\mathrm{Cay}(S_n,S))=\sum_{\sigma\in S} \tilde{\chi}_{(n-1,1)}(\sigma).
$$
One can also define whether a weighted Cayley graph on any finite group has the Aldous property with respect to a representation of the group (see \cite{C1, C3, K, PP}). In \cite{C1}, Cesi proved that the pancake graph $P_n := \mathrm{Cay}(S_n, \{r_{1,j}~|~ 2\leq j\leq n\})$ has the Aldous property, where $r_{1,j}\in S_n$ is the permutation which maps $1, 2, \ldots, j-1, j$ to $j, j-1, \ldots, 2, 1$ respectively and fixes all other points in $[n]$. Pancake graphs form a family of non-normal Cayley graphs on $S_n$ which have the Aldous property, along with the ones in the original setting of Aldous' conjecture. In \cite{PP}, Parzanchevski and Puder studied the strictly second largest eigenvalue of $\mathrm{Cay}(S_n, S)$ in the case when $S$ is a single conjugacy class of $S_n$. In \cite{HH2}, Huang, Huang and Cioab\u{a} proved that a majority of the connected normal Cayley graphs on $S_n$ ($n\ge 7$) with connection sets consisting of permutations moving at most five points possess the Aldous property. In \cite{C3}, Cesi proved that certain Cayley graphs on the Weyl group $W(B_n)$ has the Aldous property, and Kassabov proved in \cite{K} that any Cayley graph on a finite Coxeter group with respect to a specific Coxeter generating set has the Aldous property. In addition, several researchers with background in probability theory have also tried to generalize Aldous' conjecture from the probabilistic framework (see \cite[p.78]{ACDHJP}, \cite[Conjecture 1.7]{BC}, \cite[p.301]{C2} and \cite{HP}). 

In this paper, we present more classes of normal Cayley graphs on $S_n$ that have the Aldous property. Recall that the \emph{support} of a permutation $\sigma \in S_n$ is defined as $\mathrm{supp}(\sigma)=\{i\in [n]~|~\sigma(i)\ne i\}.$ For $\emptyset \ne I \subseteq \{2,3,\ldots,n-1, n\}$ and $2\leq k\leq n$, set
$$
T(n,I)=\{\sigma\in S_n \mid |\mathrm{supp}(\sigma)|\in I\}
$$ 
and
$$
T(n,k)=\{\sigma \in S_n~|~2\leq |\mathrm{supp}(\sigma)|\leq k\}.
$$
Our main results are as follows. 

\begin{thm}\label{thm3.1}
There exists a positive integer $N$ such that for every $n\ge N$ and any conjugacy class $S$ of $S_n$, the normal Cayley graph $\mathrm{Cay}(S_n, S)$ has the Aldous property if and only if $2 \le |\mathrm{supp}(\sigma)| \le n-2$ for some (and hence all) $\sigma \in S$.
\end{thm}
 
\begin{thm}\label{thm:result1}
There exists a positive integer $N$ such that for every $n\ge N$ and any $\emptyset \ne I \subset \{2,3,\ldots,n-1, n\}$ with $|I \cap \{n-1, n\}| \ne 1$, the normal Cayley graph $\mathrm{Cay}(S_n,T(n,I))$ has the Aldous property if and only if $I \cap \{n-1, n\} = \emptyset$.
\end{thm}

\begin{thm}\label{thm:result2}
There exists a positive integer $N$ such that for every $n\ge N$ and any $2\leq k\leq n$, the connected normal Cayley graph $\mathrm{Cay}(S_n,T(n,k))$ has the Aldous property.
\end{thm} 

It is worth mentioning that the case $2\leq k\leq n-2$ in Theorem \ref{thm:result2} is covered by the sufficiency part of Theorem \ref{thm:result1}. We keep the full range $2\leq k\leq n$ there since from this form one can easily see that Theorem \ref{thm:result2} is an extension of the ``normal'' version of Aldous' spectral gap conjecture. In fact, $T(n,2)$ is the conjugacy class of all transpositions of $S_n$ and $\mathrm{Cay}(S_n,T(n,2))$ is the unique normal Cayley graph covered by Aldous' spectral gap conjecture. We have $T(n,2) \subseteq T(n,k)$ for $2\leq k\leq n$ and therefore $\mathrm{Cay}(S_n,T(n,k))$ is connected indeed. 

After the completion of the proof of Theorem \ref{thm:result2}, we realized that the fact that the normal Cayley graph $\mathrm{Cay}(S_n,T(n,n-1))$ has the Aldous property can be derived from \cite[Theorem 7.1]{R} (see Remark \ref{re:re2} for details). So Theorem \ref{thm:result2} can be regarded as a corollary of Theorem \ref{thm:result1} and \cite[Theorem 7.1]{R}, modulo the trivial case of $\mathrm{Cay}(S_n,T(n,n))$. 

In \cite[Proposition 2.4]{PP}, it was proved that, for sufficiently large $n$ and any conjugacy class $S$ of $S_n$ whose elements fix at most one point, the strictly second largest eigenvalue of $\mathrm{Cay}(S_n, S)$ is attained by one of the following eight irreducible representations of $S_n$: 
$$
\rho_{(n-1,1)},\ \rho_{(n-1,1)'},\ \rho_{(n-2,2)},\ \rho_{(n-2,2)'},\ \rho_{(n-2,1,1)'},\ \rho_{(n-3,3)},\ \rho_{(n-3,2,1)},\ \rho_{(n-4,4)}.
$$
The ``only if" part of Theorem \ref{thm3.1} implies that the standard representation $\rho_{(n-1,1)}$ can be removed from this list. Also, for sufficiently large $n$, a result proved in \cite{HH} can be obtained from Theorem \ref{thm3.1} by taking $S$ to be the conjugacy class of $3$-cycles or from Theorem \ref{thm:result1} by choosing $I = \{3\}$ (see Remark \ref{re:re1} for details).

The conditions in Theorems \ref{thm3.1}, \ref{thm:result1} and \ref{thm:result2} can be stated in terms of the number of fixed points of a permutation. In particular, for $0\leq k\leq n-2$, $T(n,\{n-k\})$ is the set of permutations of $S_n$ fixing exactly $k$ points, and the normal Cayley graph $\mathrm{Cay}(S_n,T(n,\{n-k\}))$ is exactly the \emph{$k$-point-fixing graph} $\mathcal{F}(n,k)$ studied in \cite{KLW}. In particular, $\mathcal{D}_n : = T(n,\{n\})$ is the set of \emph{derangements} on $[n]$ and $\mathcal{F}(n,0) = \mathrm{Cay}(S_n, \mathcal{D}_n)$ is widely known as the \emph{derangement graph} on $[n]$. Theorem \ref{thm:result1} together with some known results from \cite{DZ, KLW, R} implies the following corollary, which asserts that, for sufficiently large $n$, $\mathcal{F}(n,0)$ and $\mathcal{F}(n,1)$ are the only graphs among $\mathcal{F}(n,k)$ ($0\leq k\leq n-2$) without the Aldous property.

\begin{cor}\label{cor1}
There exists a positive integer $N$ such that for every $n\ge N$, the $k$-point-fixing graph $\mathcal{F}(n,k)$ has the Aldous property if and only if $2\leq k\leq n-2$.
\end{cor}

The integer $N$ in all results above is no more than an integer threshold in \cite{PP} which is believed to be as small as $17$. In fact, as will be seen in the proofs of Theorems \ref{thm3.1}--\ref{thm:result2} and Corollary \ref{cor1}, the integer $N$ in these results is given by $\max\{N_0, 6\}$, $\max\{N_0, 7\}$, $\max\{N_0, 7\}$ and $\max\{N_0, 6\}$, respectively, where $N_0$ (see Lemma \ref{lem2.2}) is the integer $N_1$ in \cite[Proposition 2.3]{PP}. The larger one between this $N_1$ and the integer $N_2$ in \cite[Proposition 2.4]{PP} gives the integer threshold $N_0$ in \cite[Theorem 1.7]{PP}, and it is conjectured in \cite[Conjecture 1.8]{PP} that the smallest value of this threshold $N_0$ is $17$. Therefore, if Conjecture 1.8 in \cite{PP} is true, then in all our results above $N$ can be set to be $17$.

Theorems \ref{thm:result1} and \ref{thm:result2} and Corollary \ref{cor1} together imply that, as far as the Aldous property of $\mathrm{Cay}(S_n,T(n,I))$ for sufficiently large $n$ is concerned, where $\emptyset \ne I \subseteq \{2,3,\ldots,n-1, n\}$, the only unsettled case is the one in which $|I| \ge 2$, $|I \cap \{n-1, n\}| = 1$ and $I \neq \{2,3,\ldots,n-1\}$. Our attempts to this case suggest that the subcases where $I \cap \{n-1, n\} = \{n-1\}$ and $I \cap \{n-1, n\} = \{n\}$, respectively, may need separate treatments as they behave differently. So we propose the following two problems separately.

\begin{prob}
Give a necessary and sufficient condition for $\mathrm{Cay}(S_n,T(n,I))$ with $\{n-1\} \subset I \subset \{2,3,\ldots,n-2,n-1\}$ to have the Aldous property for sufficiently large $n$.
\end{prob}

\begin{prob}
Give a necessary and sufficient condition for $\mathrm{Cay}(S_n,T(n,I))$ with $\{n\} \subset I \subseteq \{2,3,\ldots,n-2,n\}$ to have the Aldous property for sufficiently large $n$.
\end{prob}
 
A major tool for our proofs of Theorems \ref{thm3.1} and \ref{thm:result1} is Proposition 2.3 in \cite{PP} (see Lemma \ref{lem2.2}). In particular, this proposition implies that, for sufficiently large $n$, the second largest eigenvalue of a connected normal Cayley graph as in Theorem \ref{thm:result1} satisfying $I \cap \{n-1, n\} = \emptyset$ can be obtained by comparing the two eigenvalues that correspond to the partitions $(n-1,1)$ and $(1^n)$ of $n$. As mentioned above, to prove Theorem \ref{thm:result2} we only need to consider the case when $k=n-1$, and in this case the proof can be reduced to determining the smallest eigenvalue of the complement of $\mathrm{Cay}(S_n, T(n,n-1))$. We will achieve this with the help of an asymptotic upper bound \cite{LS} on the irreducible characters of $S_n$ (see Lemma \ref{lem2.4}) and a lower bound \cite{PP} on the dimensions of the irreducible representations of $S_n$ (see Lemma \ref{lem2.5}).   

The remainder of this paper is structured as follows. In the next section, we give some basic definitions and present several known results that will play a key role in the proofs of our main results. After these preparations, the proofs of Theorems \ref{thm3.1}, \ref{thm:result1} and \ref{thm:result2} will be given in Sections \ref{sec:single}, \ref{sec:result1} and \ref{sec:result2}, respectively. The proof of Corollary \ref{cor1} will be given in Section \ref{sec:result2} as well.

\section{Preliminaries}\label{sec:pre}

All definitions in this section can be found in \cite{J,JK,Sagan}. A \emph{partition} of a positive integer $n$ is a sequence of positive integers $\gamma=(\gamma_1,\gamma_2,\ldots,\gamma_m)$ satisfying $\gamma_1\ge \gamma_2\ge\cdots\ge\gamma_m$ and $n=\gamma_1+\gamma_2+\cdots+\gamma_m$. We use $\gamma \vdash n$ to indicate that $\gamma$ is a partition of $n$ and use $c_i(\gamma)$ to denote the number of terms in $\gamma$ which are equal to $i$. A \emph{Young~diagram} is a finite collection of blocks arranged in left-justified rows, with the row sizes weakly decreasing. The Young diagram associated to the partition $\gamma=(\gamma_1,\gamma_2,\ldots,\gamma_m)$ is the one that has $m$ rows and $\gamma_i$ blocks on the $i$-th row. For instance, the following are the Young diagrams corresponding to all partitions of 4.
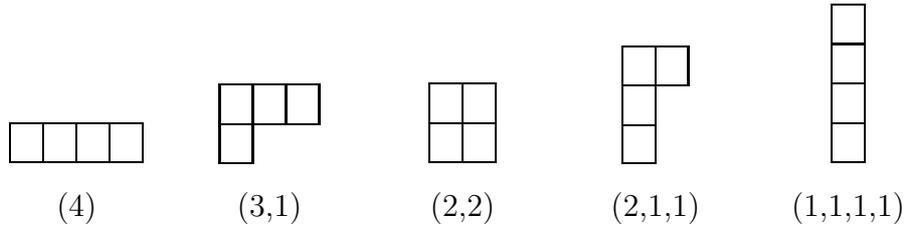
\begin{figure}[H] 
\centering
  \begin{minipage}[b]{0.15\linewidth}
    \centering
      \begin{tabular}{|l|l|l|l|}
        \hline
         & & & \\ \hline
       \end{tabular}
    \captionof*{table}{(4)}
  \end{minipage}
  \begin{minipage}[b]{0.15\linewidth}
    \centering
      \begin{tabular}{|l|l l }
        \hline
        & \multicolumn{1}{l|}{} & \multicolumn{1}{l|}{} \\ \hline
        &  &  \\ \cline{1-1}
      \end{tabular}
    \captionof*{table}{(3,1)}
  \end{minipage}
  \begin{minipage}[b]{0.15\linewidth}
    \centering
      \begin{tabular}{|l|l|}
        \hline
        &   \\ \hline
        &   \\ \hline
      \end{tabular}
    \captionof*{table}{(2,2)}
  \end{minipage}
  \begin{minipage}[b]{0.15\linewidth}
    \centering
      \begin{tabular}{|l|l}
        \hline
        & \multicolumn{1}{l|}{}  \\ \hline
        &                        \\ \cline{1-1}
        &                        \\ \cline{1-1}
      \end{tabular}
    \captionof*{table}{(2,1,1)}
  \end{minipage}
 \begin{minipage}[b]{0.15\linewidth}
   \centering
     \begin{tabular}{|l|}
         \hline
      \\ \hline
      \\ \hline
      \\ \hline
      \\ \hline
     \end{tabular}
  \captionof*{table}{(1,1,1,1)}
 \end{minipage} 
\end{figure}
\noindent Since there is a clear one-to-one correspondence between partitions and Young diagrams, we use the two terms interchangeably, and we always use Greek letters $\gamma$ and $\zeta$ to denote them. 

  %We use $\sigma_{ij}$ to denote the transposition exchanging $i$ and $j$ in $[n]$. 
The \emph{sign} of a permutation $\sigma\in S_n$, written $\mathrm{sgn}(\sigma)$, is defined to be $1$ if $\sigma$ is even and $-1$ if $\sigma$ is odd. Every permutation $\sigma\in S_n$ has a decomposition into disjoint cycles. The \emph{cycle type} of $\sigma$ is the partition of $n$ whose parts are the lengths of the cycles in its decomposition. It is widely known that two elements of $S_n$ are conjugates if and only if they have the same cycle type. This means that the conjugacy classes of $S_n$ are characterized by the cycle types and thus correspond to the partitions of $n$. Denote by $C(S_n,\gamma)$ the conjugacy class of $S_n$ that corresponds to the partition $\gamma\vdash n$. For each $\gamma\vdash n$, we define $\mathrm{sgn}(\gamma) = 1$ if all permutations in $C(S_n,\gamma)$ are even and $\mathrm{sgn}(\gamma)=-1$ otherwise.

For each partition $\zeta\vdash n$, we use $\rho_\zeta$ to denote the Specht module in $\widehat{S_n}$ that corresponds to $\zeta$. The Specht modules $\rho_{(n)}$ and $\rho_{(1^n)}$ are called the \emph{trivial} and the \emph{sign representations} of $S_n$, respectively.  We use $\chi_\zeta(\cdot)$ and $\tilde{\chi}_\zeta(\cdot)$ to denote the character and normalized character of $\rho_\zeta$, respectively. It is well known that $\chi_{(n)}(\sigma)=\tilde{\chi}_{(n)}(\sigma)= 1$ and $\chi_{(1^n)}(\sigma)=\tilde{\chi}_{(1^n)}(\sigma)=\mathrm{sgn}(\sigma)$ for any $\sigma\in S_n$. Since $\chi_{\zeta}(\cdot)$ (respectively, $\tilde{\chi}_{\zeta}(\cdot)$) is a class function on $S_n$, we use $\chi_\zeta(\gamma)$ (respectively, $\tilde{\chi}_\zeta(\gamma)$) to denote the value of $\chi_\zeta(\cdot)$ (respectively, $\tilde{\chi}_\zeta(\cdot)$) on the conjugacy class $C(S_n,\gamma)$. We always use $\iota$ to denote the identity element of $S_n$. Note that $\chi_\zeta(\iota)$ equals the dimension of $\rho_\zeta\in \widehat{S_n}$ for any $\zeta\vdash n$.

It is known \cite{F} that the character of any $\rho_\zeta\in \widehat{S_n}$ on any conjugacy class of $S_n$ is an integer with absolute value at most the dimension of $\rho_\zeta$. Hence $\tilde{\chi}_\zeta(\gamma)$ is a rational number in the interval $[-1,1]$ for all $\zeta,\gamma\vdash n$. For the convenience of the reader and in order to provide self-contained proofs, we include the following Table~\ref{tab:tab1} from \cite{PP}, which gives the dimensions and characters of some irreducible representations of $S_n$.  

\begin{table}[h]
\centering
\begin{tabular}{c c c}
%\vspace{0.1cm}
\toprule
$~~~~\zeta\vdash n~~~$ & $~~~~\dim \rho_\zeta=\chi_\zeta(\iota)~~~~$ & $~~~~\chi_\zeta(\gamma)~\text{with}~c_i(\gamma)=c_i~~~~$\\\toprule
$(n)$ & $1$ & $1$\\ \midrule
$(n-1,1)$ & $n-1$ & $c_1-1$\\\midrule
$(n-2,2)$ & $\frac{n(n-3)}{2}$ & $\frac{c_1(c_1-3)}{2}+c_2$\\\midrule
$(n-2,1,1)$& $\frac{(n-1)(n-2)}{2}$ & $\frac{(c_1-1)(c_1-2)}{2}-c_2$\\ \midrule
$(n-3,3)$& $\frac{n(n-1)(n-5)}{6}$ & $\frac{c_1(c_1-1)(c_1-5)}{6}+(c_1-1)c_2+c_3$\\ \midrule
$(n-3,2,1)$ & $\frac{n(n-2)(n-4)}{3}$ & $\frac{c_1(c_1-2)(c_1-4)}{3}-c_3$ \\ 
%$(n-4,4)$ & $\frac{n(n-1)(n-2)(n-7)}{24}$ & $\begin{aligned}
%	\frac{c_1(c_1-1)(c_1-2)(c_1-7)}{24}+\frac{(c_1^2-3c_1-1)c_2}{2}\\+\frac{c_2^2}{2}+(c_1-1)c_3+c_4~~~~~~~~~~~~~~
%\end{aligned}$ \\
\bottomrule
\end{tabular}
\vspace{0.2cm}
\caption{Dimensions and characters of some irreducible representations of $\widehat{S_n}$}
\label{tab:tab1}
\end{table}

The \emph{conjugate} or \emph{transpose} of a partition $\zeta=(\zeta_1,\zeta_2,\ldots,\zeta_m)\vdash n$ is defined as $\zeta'=(\zeta_1',~\zeta_2',~\ldots,\zeta_h') \vdash n$, where $\zeta_i'$ is the length of the $i$-th column of $\zeta$. In other words, the Young diagram of $\zeta'$ is just the transpose of that of $\zeta$. % so sometimes $\zeta'$ is written as $\zeta^t$.  
The relation between $\chi_\zeta(\cdot)$ and $\chi_{\zeta'}(\cdot)$ is reflected in the following lemma.

\begin{lem}{\upshape (\cite[2.1.8]{JK}) }\label{lem2.1}
  For any $\zeta, \gamma \vdash n$, we have 
  \begin{equation*}
	\chi_{\zeta'}(\gamma)=\mathrm{sgn}(\gamma)\cdot\chi_\zeta(\gamma).
  \end{equation*} 
\end{lem}

Let $E(n)$ and $O(n)$ be the numbers of even and odd derangements on $[n]$, respectively. The next two lemmas will be used in our proof of  Theorem \ref{thm:result1}. 

\begin{lem}{\upshape (\cite{O}) }\label{eo}
$E(n)-O(n)=(-1)^{(n-1)}(n-1).$ 
\end{lem}

\begin{lem}{\upshape (\cite[Proposition 2.3]{PP}) }\label{lem2.2}
There exists a positive integer $N_{0}$ such that for every $n\ge N_{0}$ and any $\gamma\vdash n$ with $c_1(\gamma)\ge 2$, we have $$\max_{\substack{\zeta\vdash n \\ \zeta \notin \{(n), (1^n)\}}} \tilde{\chi}_\zeta(\gamma)=\tilde{\chi}_{(n-1,1)}(\gamma).$$ 
\end{lem}

Our main tool for proving the case $k=n-1$ in Theorem \ref{thm:result2} is the following asymptotically sharp bound for the characters of $S_n$ due to Larsen and Shalev \cite{LS}.   

\begin{lem}{\upshape (\cite[Theorem 1.3]{LS}) }\label{lem2.4}
Let $\gamma$ be a partition of $n$ and let $f=\max\{c_1(\gamma), 1\}$. Then for all $\zeta \vdash n$,
$$
|\chi_\zeta(\gamma)| \leq \chi_\zeta(\iota)^{1-\frac{\log (n / f)}{2 \log n}+ \varepsilon_n}=\chi_\zeta(\iota)^{\frac{1}{2}+\frac{\log f}{2 \log n}+\varepsilon_n},
$$ where $\varepsilon_n$ is a real number tending to $0$ as $n\to \infty$.
\end{lem}

Another key ingredient for proving the case $k=n-1$ in Theorem \ref{thm:result2} is the following estimation of the dimensions of irreducible representations of $S_n$.

\begin{lem}{\upshape (\cite[Lemma 2.6]{PP}) }\label{lem2.5}
Let $n\ge 13$ and $\zeta\vdash n$ correspond to a Young diagram with at least three blocks outside the first row and at least three blocks outside the first column. Then the dimension $\chi_\zeta(\iota)$ of $\rho_\zeta\in \widehat{S_n}$ satisfies 
$$
\chi_\zeta(\iota)\ge n^{2.05}.
$$ 
\end{lem}

\section{Proof of Theorem \ref{thm3.1}}
\label{sec:single}

%\begin{equation}\label{pp}
%	\max_{\substack{\zeta\vdash n \\ \zeta\ne (n)~\text{or}~(1^n)}}\tilde{\chi}_\zeta(\gamma)=\max_{\zeta\in \mathfrak{eight}_n}\tilde{\chi}_\zeta(\gamma),
%\end{equation} 

Note that $C(S_n,\gamma) \ne \{\iota\}$ if and only if $\gamma \ne (1^n)$, which in turn is true if and only if $|\mathrm{supp}(\sigma)| \ge 2$ for any $\sigma \in C(S_n,\gamma)$. Note also that $c_1(\gamma)\ge 2$ if and only if $|\mathrm{supp}(\sigma)| \le n-2$ for any $\sigma \in C(S_n,\gamma)$. So Theorem \ref{thm3.1} can be restated as follows: There is a positive integer $N$ such that for every $n\ge N$ and any $(1^n)\ne \gamma\vdash n$, the Cayley graph $\mathrm{Cay}(S_n, C(S_n,\gamma))$ has the Aldous property if and only if $c_1(\gamma)\ge 2$. We prove this statement in the rest of this section.
 
\begin{proof} 
Since $\mathrm{Cay}(S_n, C(S_n,\gamma))$ is normal, by Proposition \ref{prop1.2} its eigenvalue corresponding to the partition $\zeta\vdash n$ is 
\[ \lambda_\zeta:= \sum_{\sigma\in C(S_n,\gamma)}\tilde{\chi}_\zeta(\sigma)=|C(S_n,\gamma)|\cdot \tilde{\chi}_\zeta(\gamma).\]
 Note that for any partition $\gamma \ne (1^n)$ of $n$, the subgroup $\langle C(S_n,\gamma)\rangle$ is either $S_n$ or $A_n$. Assume first $\langle C(S_n,\gamma)\rangle=S_n$. Then every permutation in $C(S_n,\gamma)$ is odd, and $|C(S_n,\gamma)|$ is a simple eigenvalue of $\mathrm{Cay}(S_n, C(S_n,\gamma))$, which is attained by the trivial representation $\rho_{(n)}$ as $\lambda_{(n)}=|C(S_n,\gamma)|\cdot \tilde{\chi}_{(n)}(\gamma)=|C(S_n,\gamma)|$. The eigenvalue corresponding to the sign representation $\rho_{(1^n)}$ is $\lambda_{(1^n)}=|C(S_n,\gamma)|\cdot \tilde{\chi}_{(1^n)}(\gamma)=-|C(S_n,\gamma)|$. According to the Perron-Frobenius Theorem \cite[Theorem 8.4.4]{HJ1}, we know that this is the smallest eigenvalue of $\mathrm{Cay}(S_n, C(S_n,\gamma))$. Now assume $\langle C(S_n,\gamma)\rangle=A_n$. Then the largest eigenvalue $|C(S_n,\gamma)|$ of $\mathrm{Cay}(S_n, C(S_n,\gamma))$ has multiplicity 2 and is attained simultaneously by the trivial representation $\rho_{(n)}$ and the sign representation $\rho_{(1^n)}$. Therefore, for any partition $\gamma \ne (1^n)$ of $n$, the strictly second largest eigenvalue of $\mathrm{Cay}(S_n,C(S_n,\gamma))$ can be expressed as 
\begin{equation}\label{thm3.1 eq}
\max_{\substack{\zeta\vdash n\\\zeta \notin \{(n), (1^n)\}}}\lambda_\zeta = |C(S_n,\gamma)|\max_{\substack{\zeta\vdash n\\\zeta \notin \{(n), (1^n)\}}}\tilde{\chi}_\zeta(\gamma).
 	\end{equation}

If $c_1(\gamma)\ge 2$, then by Lemma \ref{lem2.2} there is a positive integer $N_{0}$ such that for every $n\ge N_{0}$, the maximum on the right-hand side of (\ref{thm3.1 eq}) is attained by $\zeta=(n-1,1)$. In other words, if $c_1(\gamma)\ge 2$, then $\mathrm{Cay}(S_n, C(S_n,\gamma))$ has the Aldous property for every $n\ge N_{0}$. 

In the remaining proof we assume that $n\ge 6$ and $(1^n)\ne \gamma\vdash n$ is such that $c_1(\gamma)\leq 1$. We aim to prove that $\mathrm{Cay}(S_n, C(S_n,\gamma))$ does not have the Aldous property in this case.

If $c_1(\gamma)=1$, then a direct computation using Table \ref{tab:tab1} leads to $\tilde{\chi}_{(n-1,1)}(\gamma)=0$ and thus $\lambda_{(n-1,1)}=0$, while the value of (\ref{thm3.1 eq}) is at least  
\begin{align*} 
\max\{\lambda_{(n-3,3)},~\lambda_{(n-3,2,1)}\} & = |C(S_n,\gamma)|\cdot \max\{\tilde{\chi}_{(n-3,3)}(\gamma),~\tilde{\chi}_{(n-3,2,1)}(\gamma)\}\\
& = |C(S_n,\gamma)|\cdot \max\left\{ \frac{6c_3(\gamma)}{n(n-1)(n-5)},~\frac{3-3c_3(\gamma)}{n(n-2)(n-4)}\right\},
\end{align*}
which is positive as $c_3(\gamma)\ge 0$. This means that the maximum on the left-hand side of (\ref{thm3.1 eq}) is not attained by $\zeta=(n-1,1)$, and thus $\mathrm{Cay}(S_n, C(S_n,\gamma))$ with $c_1(\gamma)=1$ does not possess the Aldous property for every $n\ge 6$.

Finally, let us consider the case when $c_1(\gamma)=0$. In this case, we have \(\tilde{\chi}_{(n-1,1)}(\gamma)=-\frac{1}{n-1}<0\) and hence $\lambda_{(n-1,1)}<0$. If $\mathrm{sgn}(\gamma)=-1$, then according to Lemma \ref{lem2.1} the value of (\ref{thm3.1 eq}) is at least 
$$
\lambda_{(n-1,1)'}=|C(S_n,\gamma)|\cdot \tilde{\chi}_{(n-1,1)'}(\gamma)=|C(S_n,\gamma)|\cdot\mathrm{sgn}(\gamma)\tilde{\chi}_{(n-1,1)}(\gamma)>0.
$$ 
If $\mathrm{sgn}(\gamma)=1$, then the value of (\ref{thm3.1 eq}) is at least 
\begin{align*}
\max\{\lambda_{(n-2,2)},~\lambda_{(n-2,1,1)'}\}
& = |C(S_n,\gamma)|\cdot\max \left\{\tilde{\chi}_{(n-2,2)}(\gamma),~\tilde{\chi}_{(n-2,1,1)'}(\gamma)\right\}\\
& = |C(S_n,\gamma)|\cdot\max \left\{ \frac{2c_2(\gamma)}{n(n-3)},~ \frac{2-2c_2(\gamma)}{(n-1)(n-2)}\right\},	
\end{align*} 
which is positive as $c_2(\gamma)\ge 0$. Thus the maximum on the left-hand side of (\ref{thm3.1 eq}) is not attained by $\zeta=(n-1,1)$. In other words, if $c_1(\gamma)=0$, then $\mathrm{Cay}(S_n, C(S_n,\gamma))$ does not possess the Aldous property for every $n\ge 6$.

To sum up, we have proved that for $n\ge N := \max\{N_{0},6\}$, the normal Cayley graph $\mathrm{Cay}(S_n, C(S_n,\gamma))$ with $(1^n)\ne\gamma\vdash n$ has the Aldous property if and only if $c_1(\gamma)\ge 2$.
\end{proof}

\section{Proof of Theorem \ref{thm:result1}}
\label{sec:result1}

\begin{proof}
We prove the sufficiency and necessity separately.

\smallskip
\textsf{Sufficiency.}  Suppose that $I \cap \{n-1, n\} = \emptyset$. Then $\emptyset \ne I \subseteq \{2,3,\ldots,n-2\}$. By the definition of $T(n,I)$, the subgroup $\langle T(n,I)\rangle$ is either $S_n$ or $A_n$. Since $\mathrm{Cay}(S_n,T(n,I))$ is normal, by Proposition \ref{prop1.2} its eigenvalue corresponding to the partition $\zeta\vdash n$ is 
$$
\lambda_\zeta:=\sum_{\sigma\in T(n,I)}\tilde{\chi}_\zeta(\sigma)=\sum_{\substack{ \gamma\vdash n\\ n-c_1(\gamma)\in I }}  |C(S_n,\gamma)|\cdot \tilde{\chi}_\zeta(\gamma).
$$
 
\smallskip
\textsf{Case 1.} $\langle T(n,I)\rangle=A_n$.
\smallskip

In this case we have $I=\{3\}$ and $T(n,I)=C\left(S_n, (3,1^{n-3})\right)$ is the single conjugacy class of all $3$-cycles in $S_n$. By Theorem \ref{thm3.1} and its proof, there is a positive integer $N_1 := \max\{N_{0},6\}$ such that for every $n\ge N_1$, the Cayley graph $\mathrm{Cay}(S_n, T(n,\{3\}))$ has the Aldous property, where $N_0$ is the positive integer given in Lemma \ref{lem2.2}.

%Similarly to  The eigenvalue $|T(n,I)|$ has multiplicity $2$ and is attained simultaneously by $\rho_{(n)}$ and $\rho_{(1^n)}$ . Since $\mathrm{Cay}(A_n,T(n,I))$ is connected, $|T(n,I)|$ is an eigenvalue of $\mathrm{Cay}(S_n, T(n,I))$ with multiplicity $2$ and all other eigenvalues of $\mathrm{Cay}(S_n, T(n,I))$ are strictly smaller than $|T(n,I)|$.  Thus the largest eigenvalue of $\mathrm{Cay}(S_n,T(n,I))$ except $|T(n,I)|$ equals \begin{equation}\label{eq1.5}
%	\big|C\big(S_n, (3,1^{n-3})\big)\big|\max_{\substack{\zeta\vdash n\\ \zeta\ne (n) ~\text{or}~(1^n)}}   \tilde{\chi}_\zeta\big((3,1^{n-3})\big).
%\end{equation}
%By Lemma \ref{lem2.2}, there is a positive integer $N_1$ such that for every $n\ge N_1$, $$\max_{\substack{\zeta\vdash n\\ \zeta\ne (n) ~\text{or}~(1^n)}} \tilde{\chi}_\zeta\big((3,1^{n-3})\big)=\tilde{\chi}_{(n-1,1)}\big((3,1^{n-3})\big)>0.$$  Therefore, for every $n\ge N_1$, the maximum in (\ref{eq1.5}) is achieved by $\zeta=(n-1,1)$, which means the largest eigenvalue of $\mathrm{Cay}(S_n,T(n,I))$ except $|T(n,I)|$ is achieved by the standard representation of $S_n$. 

\smallskip
\textsf{Case 2.} $\langle T(n,I)\rangle=S_n$.
\smallskip

It is clear that the eigenvalue of $\mathrm{Cay}(S_n, T(n,I))$ corresponding to the trivial representation $\rho_{(n)}$ is $\lambda_{(n)}=|T(n,I)|$, which is simple as $\mathrm{Cay}(S_n, T(n,I))$ is connected. The second largest eigenvalue of $\mathrm{Cay}(S_n, T(n,I))$, which is also the strictly second largest eigenvalue of $\mathrm{Cay}(S_n, T(n,I))$, is given by
\begin{equation}\label{eq1.6}
	\max_{\substack{\zeta\vdash n\\ \zeta\ne (n)}} \lambda_\zeta=\max_{\substack{\zeta\vdash n\\ \zeta\ne (n)}}\sum_{\substack{ \gamma\vdash n\\ n-c_1(\gamma)\in I }}  |C(S_n,\gamma)|\cdot \tilde{\chi}_\zeta(\gamma).
\end{equation}
To complete the proof, it suffices to prove that the maximum in (\ref{eq1.6}) is attained by $\zeta=(n-1,1)$.

Since $I$ is a nonempty subset of $\{2,3,\ldots,n-2\}$, every partition $\gamma\vdash n$ with $n-c_1(\gamma)\in I$ satisfies $c_1(\gamma)\ge 2$. Thus, by Lemma \ref{lem2.2}, there is a positive integer $N_0$ such that for every $n\ge N_{0}$ and any $\gamma\vdash n$ with $n-c_1(\gamma)\in I$, we have
$$
\max_{\substack{\zeta\vdash n\\\zeta \notin \{(n), (1^n)\}}}\tilde{\chi}_\zeta(\gamma)=\tilde{\chi}_{(n-1,1)}(\gamma).
$$ 
This implies that, for every $n\ge N_{0}$,
\begin{align*}
	\max_{\substack{\zeta\vdash n\\ \zeta \notin \{(n), (1^n)\}}} \lambda_\zeta\ & =  \max_{\substack{\zeta\vdash n\\ \zeta\notin \{(n), (1^n)\}}} \sum_{\substack{\gamma\vdash n\\ n-c_1(\gamma)\in I}} |C(S_n,\gamma)|\cdot \tilde{\chi}_{\zeta}(\gamma) \\
	& = \sum_{\substack{\gamma\vdash n\\ n-c_1(\gamma)\in I}} |C(S_n,\gamma)|\cdot \tilde{\chi}_{(n-1,1)}(\gamma) \\
	& = \lambda_{(n-1,1)}. 
\end{align*}
Thus, for every $n\ge N_{0}$, the second largest eigenvalue of $\mathrm{Cay}(S_n, T(n,I))$ given in formula (\ref{eq1.6}) equals
\begin{equation}\label{eq1.7}
	\max\{\lambda_{(n-1,1)},~\lambda_{(1^n)}\}=\max_{\zeta \in \{(n-1,1), (1^n)\}} \sum_{\substack{\gamma\vdash n\\ n-c_1(\gamma)\in I}} |C(S_n,\gamma)|\cdot \tilde{\chi}_{\zeta}(\gamma).
\end{equation}
In the following we prove that the eigenvalue $\lambda_{(n-1,1)}$ is always greater than or equal to $\lambda_{(1^n)}$ for any $n\ge 7$ and $\emptyset \ne I\subseteq \{2,3,\ldots,n-2\}$ with $I \ne \{3\}$. This will be accomplished by giving explicit expressions for $\lambda_{(n-1,1)}$ and $\lambda_{(1^n)}$ and then comparing them on a term-by-term basis.

Let us first calculate the values of $\lambda_{(n-1,1)}$ and $\lambda_{(1^n)}$. When $\zeta=(n-1,1)$, we have $\chi_\zeta(\gamma)=c_1(\gamma)-1$ and $\chi_{\zeta}(\iota)=n-1$ by Table \ref{tab:tab1}. So 
\begin{align}
\lambda_{(n-1,1)} & =\ \sum_{\substack{\gamma\vdash n\\ n-c_1(\gamma)\in I}} |C(S_n,\gamma)|\cdot \frac{c_1(\gamma)-1}{\chi_\zeta(\iota)} \nonumber\\
& =\ \sum_{t\in I}\frac{n-t-1}{n-1}\cdot \sum_{\substack{\gamma\vdash n\\ n-c_1(\gamma)=t}} |C(S_n,\gamma)| \nonumber \\
& =\ \sum_{t\in I} \frac{n-t-1}{n-1}\cdot \binom{n}{t} \cdot D(t) \label{eq1.8} \\
& >\ 0, \label{eq18}
\end{align}
where $D(t)=t! \cdot \sum_{i=0}^{t} \frac{(-1)^i}{i!} $ is the number of derangements on $[t]$. Note that $D(t)\ge \frac{t!}{3}$ for every integer $t\ge 3$. On the other hand, if $\zeta=(1^n)$, then $\chi_\zeta(\gamma)=\mathrm{sgn}(\gamma)$ and $\chi_\zeta(\iota)=1$. Thus, 
\begin{align}
\lambda_{(1^n)} &=\ \sum_{\substack{\gamma\vdash n\\ n-c_1(\gamma)\in I}} |C(S_n,\gamma)|\cdot \mathrm{sgn}(\gamma)\nonumber\\
& =\ \sum_{t\in I} \sum_{\substack{\gamma\vdash n\\ n-c_1(\gamma)=t}} |C(S_n,\gamma)|\cdot \mathrm{sgn}(\gamma) \nonumber \\
& =\ \sum_{t\in I} \binom{n}{t} \left(E(t)-O(t)\right) \nonumber \\
& =\ \sum_{t\in I} (-1)^{t-1}(t-1)\binom{n}{t}, \label{eq1.10}
\end{align}
where $E(t)$ and $O(t)$ are the numbers of even and odd derangements on $[t]$, respectively, and the last step follows from Lemma \ref{eo}. 

Now we prove that when $n\ge 7$, for any $\emptyset \ne I \subseteq \{2,3,\ldots,n-2\}$ with $I \ne \{3\}$, the value of formula (\ref{eq1.8}) is no less than that of formula (\ref{eq1.10}). Define 
$$
A(t):=\frac{n-t-1}{n-1} \cdot \binom{n}{t} \cdot D(t),\quad B(t):=(-1)^{t-1}(t-1)\binom{n}{t}
$$ 
for $2\leq t\leq n-2$.
By straightforward computations, one can verify that $A(2)-B(2)=n(n-2)$, $A(3)-B(3)=-n(n-2)$, and $A(4)-B(4)=\frac{1}{2}n(n-2)(n-3)(n-4)>n(n-2)$ for $n\ge 7$. For $5\leq t\leq n-2$, we have
\begin{align*}
A(t)-B(t) & = \frac{n-t-1}{n-1}\cdot \binom{n}{t} \cdot D(t)-(-1)^{t-1}(t-1)\binom{n}{t} \\
& \ge  \frac{n-t-1}{n-1}\cdot \binom{n}{t} \cdot D(t)-(t-1)\binom{n}{t} \\
& >  \frac{n-t-1}{n-1}\cdot \binom{n}{t} \cdot \frac{t!}{3}-(t-1)\binom{n}{t} \\
& =  \frac{n-t-1}{n-1}\cdot \binom{n}{t} \cdot \frac{t(t-1)(t-2)!}{3}-(t-1)\binom{n}{t} \\
& \ge  \frac{n-t-1}{n-1}\cdot 2t(t-1)\cdot \binom{n}{t} -(t-1)\binom{n}{t} \\
& =  \frac{-2t^2+2nt-2t-n+1}{n-1}(t-1)\binom{n}{t},
\end{align*}
where in the second last step we used the fact that $(t-2)!\ge 6$ when $t\ge 5$. Let $f(t)=-2t^2+2nt-2t-n+1$. It can be verified that $f(t)>0$ for every positive integer $t$ in the interval $[5,n-2]$ and that the minimum value of $f(t)$ in this interval is $f(n-2)=n-3$. So, continuing the estimation above, we obtain
\begin{align*}
A(t)-B(t) & >  \frac{n-3}{n-1}(t-1)\binom{n}{t} \\
	& = \frac{n(n-3)(t-1)}{2}\cdot \frac{(n-2)(n-3)\cdots(n-t+1)}{t(t-1)\cdots3} \\
	& \ge  \frac{n(n-3)(t-1)}{2} \\
	& \ge  2n(n-3) \\
	& >  n(n-2).
\end{align*}

So far we have proved that if $n\ge 7$, then $A(3)-B(3)=- n(n-2)$ and $A(t)-B(t)\ge n(n-2)$ for any integer $t$ with $3\ne t\in [2,n-2]$. Therefore, when $n\ge 7$, for any $\emptyset \ne I\subseteq \{2,3,\ldots, n-2\}$ with $I \ne \{3\}$, the value of formula (\ref{eq1.8}) is greater than or equal to that of formula (\ref{eq1.10}). This implies that the maximum in formula (\ref{eq1.7}) is always attained by $\zeta=(n-1,1)$ when $n\ge 7$, and thus the maximum in (\ref{eq1.6}) is achieved by $\zeta=(n-1,1)$ for every $n \ge N_{2} := \max\{N_{0},7\}$.

To sum up, we have proved that for every $n\ge N := \max\{N_{1}, N_{2}\} = \max\{N_{0}, 7\}$ and any $\emptyset \ne I\subseteq \{2,3,\ldots, n-2\}$, the normal Cayley graph $\mathrm{Cay}(S_n,T(n,I))$ has the Aldous property, where $N_0$ is the positive integer given in Lemma \ref{lem2.2}. This establishes the sufficiency part of Theorem \ref{thm:result1}.

\smallskip
\textsf{Necessity.} Since $|I \cap \{n-1,n\}| \ne 1$ by our assumption, to establish the necessity it suffices to prove the following statement: There exists a positive integer $N$ such that for every $n\ge N$ and any $\{n-1,n\}\subseteq J \subset \{2,3,\ldots,n-1,n\}$, $\mathrm{Cay}(S_n,T(n,J))$ does not have the Aldous property.

Now suppose that $\{n-1,n\}\subseteq J \subset \{2,3,\ldots,n-1,n\}$. Set $I=\{2,3,\ldots,n-1,n\}\setminus J$. Then $\emptyset\ne I\subseteq\{2,3,\ldots,n-2\}$ and $\{T(n,J), T(n,I)\}$ is a partition of $S_n\setminus \{\iota\}$, where as before $\iota$ denotes the identity of $S_n$. Hence $\mathrm{Cay}(S_n,T(n,J))$ and $\mathrm{Cay}(S_n,T(n,I))$ are complements of each other. It is clear that $\mathrm{Cay}(S_n,T(n,J))$ is connected and its largest eigenvalue $|T(n,J)|=n!-|T(n,I)|-1$ is achieved by the trivial representation $\rho_{(n)}$ of $S_n$. According to Proposition \ref{prop1.2}, any other eigenvalue of $\mathrm{Cay}(S_n,T(n,J))$ can be expressed as 
$
\sum_{\sigma\in T(n,J)} \tilde{\chi}_{\zeta}(\sigma)
$
for some $(n)\ne\zeta\vdash n$. Since the complete graph $\mathrm{Cay}(S_n, S_n\setminus \{\iota\})$ has eigenvalues $n! - 1$ with multiplicity $1$ and $-1$ with multiplicity $n! - 1$, by Proposition \ref{prop1.2}, we have $\sum_{\sigma\in S_n\setminus \{\iota\}} \tilde{\chi}_{\zeta}(\sigma) = -1$ for any $(n)\ne\zeta\vdash n$. The fact that $\{T(n,J), T(n,I)\}$ is a partition of $S_n\setminus \{\iota\}$ enables us to write this as  
\begin{equation}\label{eq:cor6.2}
\sum_{\sigma\in T(n,J)} \tilde{\chi}_{\zeta}(\sigma)+\sum_{\sigma\in T(n,I)} \tilde{\chi}_\zeta(\sigma)=-1
\end{equation}
for any $(n)\ne\zeta\vdash n$. Note that $\sum_{\sigma\in T(n,I)} \tilde{\chi}_\zeta(\sigma)$ is an eigenvalue of $\mathrm{Cay}(S_n,T(n,I))$.

\smallskip
\textsf{Case 3.} $I\ne \{3\}$.
\smallskip

In this case, $\mathrm{Cay}(S_n,T(n,I))$ is connected and by the sufficiency proved above, there exists a positive integer $N$ such that
\begin{equation}\label{eq:cor6.3}
	\max_{\substack{\zeta\vdash n \\ \zeta\ne (n)}} \sum_{\sigma\in T(n,I)} \tilde{\chi}_{\zeta}(\sigma)=\sum_{\sigma\in T(n,I)} \tilde{\chi}_{(n-1,1)}(\sigma)
\end{equation}
whenever $n\ge N$. Since the sum of the eigenvalues of $\mathrm{Cay}(S_n,T(n,I))$ is $0$ and by \eqref{eq18} the second largest eigenvalue of $\mathrm{Cay}(S_n,T(n,I))$ as shown on the right-hand side of \eqref{eq:cor6.3} is positive, it follows that $\mathrm{Cay}(S_n,T(n,I))$ has at least three distinct eigenvalues and at least one of them is negative. This together with (\ref{eq:cor6.2}) and (\ref{eq:cor6.3}) implies that $\mathrm{Cay}(S_n,T(n,J))$ has at least three distinct eigenvalues and the smallest one of them is attained by the standard representation of $S_n$.   %which equals $$\sum_{\sigma\in T(n,J)} \tilde{\chi}_{(n-1,1)}(\sigma)=-1-\sum_{\sigma\in T(n,I)} \tilde{\chi}_{(n-1,1)}(\sigma).$$ 
Hence $\mathrm{Cay}(S_n,T(n,J))$ does not have the Aldous property when $n\ge N$.

\smallskip
\textsf{Case 4.} $I=\{3\}$.
\smallskip

In this case, $\mathrm{Cay}(S_n,T(n,I))$ is the union of two copies of $\mathrm{Cay}(A_n,T(n,I))$. So the largest eigenvalue $|T(n,I)|$ of $\mathrm{Cay}(S_n,T(n,I))$ has multiplicity $2$ and is attained simultaneously by the trivial and sign representations of $S_n$. Thus, by (\ref{eq:cor6.2}), we know that the sign representation attains the smallest eigenvalue of $\mathrm{Cay}(S_n,T(n,J))$, which is $-1-|T(n,I)|$. Moreover, by the sufficiency proved above, whenever $n\ge N$ the strictly second largest eigenvalue of $\mathrm{Cay}(S_n,T(n,I))$ is
\begin{equation}\label{eq:cor6.4}
	\max_{\substack{\zeta\vdash n \\ \zeta \notin \{(n),(1^n)\}}} \sum_{\sigma\in T(n,I)} \tilde{\chi}_{\zeta}(\sigma)=\sum_{\sigma\in T(n,I)} \tilde{\chi}_{(n-1,1)}(\sigma)>0.
\end{equation}
Similarly to the case above, we obtain that $\mathrm{Cay}(S_n,T(n,I))$ has at least one negative eigenvalue, say $\lambda$, as the sum of its eigenvalues is $0$. This together with (\ref{eq:cor6.2}) and (\ref{eq:cor6.4}) implies that the second smallest eigenvalue of $\mathrm{Cay}(S_n,T(n,J) )$  is attained by the standard representation of $S_n$, which is larger than $-1-|T(n,I)|$ but strictly smaller than $-1-\lambda$. Thus, the strictly second largest eigenvalue of $\mathrm{Cay}(S_n,T(n,J))$ is not attained by the standard representation of $S_n$. In other words, $\mathrm{Cay}(S_n,T(n,J))$ does not have the Aldous property when $n\ge N$.
\end{proof}

\begin{re}
\label{re:re1}
In \cite[Theorem 3.4]{HH}, Huang and Huang determined the second largest eigenvalue of the \emph{complete alternating group graph} $\mathrm{Cay}(A_n, \{(i~ j~ k), (i~ k~ j)~|~1\leq i<  j<k\leq n\})$. Note that this graph is exactly $\mathrm{Cay}(A_n,T(n,I))$ with $I=\{3\}$. Since $\mathrm{Cay}(S_n,T(n,\{3\}))$ is the union of two copies of $\mathrm{Cay}(A_n,T(n,\{3\}))$, it has the same eigenvalues as the latter but with the multiplicity of each eigenvalue doubled. Thus, for sufficiently large $n$, we can obtain \cite[Theorem 3.4]{HH} from the sufficiency part of Theorem \ref{thm:result1} by choosing $I = \{3\}$ or from Theorem \ref{thm3.1} by taking $S$ to be the conjugacy class of $3$-cycles. In fact, by Theorem \ref{thm3.1} or \ref{thm:result1}, $\mathrm{Cay}(S_n,T(n,\{3\}))$ has the Aldous property for sufficiently large $n$. Using this and Table \ref{tab:tab1}, we obtain that, for sufficiently large $n$,
\begin{align*}
\lambda_2(\mathrm{Cay}(A_n,T(n,\{3\})) & =\ \lambda_3(\mathrm{Cay}(S_n,T(n,\{3\})) \\
& =\  \sum\limits_{\sigma\in T(n,\{3\})}\tilde{\chi}_{(n-1,1)}(\sigma) \\
& =\  |T(n,\{3\})|\cdot \tilde{\chi}_{(n-1,1)}\big((3,1^{n-3})\big) \\
& =\  |T(n,\{3\})|\cdot\frac{n-4}{n-1} \\
& =\  \frac{1}{3}n(n-2)(n-4),
\end{align*}
which is exactly what is stated in \cite[Theorem 3.4]{HH}.  
\end{re}

\section{Proofs of Theorem \ref{thm:result2} and Corollary \ref{cor1}} 
\label{sec:result2}

We prove Theorem \ref{thm:result2} first.

\begin{proof}
It is clear that for any $2\leq k\leq n$ we have
\begin{equation*}
T(n,k)=\mathop{\bigcup}\limits_{\substack{\gamma \vdash n \\ 2\leq n-c_1(\gamma)\leq k}} C(S_n,\gamma)
\end{equation*}
and hence $\mathrm{Cay}(S_n,T(n,k))$ is normal.
Thus, by Proposition \ref{prop1.2}, the eigenvalues of $\mathrm{Cay}(S_n,T(n,k))$ can be expressed as 
\begin{equation}\label{eq1}
	\sum_{\sigma\in T(n,k)} \tilde{\chi}_\zeta(\sigma)=\sum_{\substack{\gamma\vdash n \\ 2\leq n- c_1(\gamma)\leq k}} |C(S_n,\gamma)| \cdot \tilde{\chi}_\zeta(\gamma),
\end{equation}
where $\zeta \vdash n$ is running over all partitions of $n$. Moreover, $\mathrm{Cay}(S_n,T(n,k))$ is connected as $T(n,2) \subseteq T(n,k)$ and $T(n,2)$ is a generating subset of $S_n$. So the largest eigenvalue $|T(n,k)|$ of $\mathrm{Cay}(S_n,T(n,k))$ is simple and is attained by the trivial representation $\rho_{(n)}$. In fact, if $\zeta=(n)$, then $\chi_\zeta(\cdot)$ is the trivial character and $\tilde{\chi}_\zeta(\cdot) = 1$. Hence  (\ref{eq1}) evaluates to $|T(n,k)|$ when $\zeta=(n)$. Thus, the second largest eigenvalue of $\mathrm{Cay}(S_n,T(n,k))$, which is also the strictly second largest eigenvalue of $\mathrm{Cay}(S_n,T(n,k))$, is equal to 
\begin{equation}\label{eq2}
	\max_{\substack{\zeta\vdash n\\ \zeta\ne (n)}} \sum_{\substack{\gamma\vdash n\\ 2\leq n-c_1(\gamma)\leq k}} |C(S_n,\gamma)|\cdot \tilde{\chi}_{\zeta}(\gamma).
\end{equation} 
To prove Theorem \ref{thm:result2}, it suffices to show that the maximum in (\ref{eq2}) is achieved by $\zeta=(n-1,1)$.  

\smallskip
\textsf{Case 1.} $k = n$, where $n \ge 4$.
\smallskip

In this case, we have $T(n,k)=S_n\setminus \{\iota\}$ and $\mathrm{Cay}(S_n,T(n,k))$ is the complete graph, which has two distinct eigenvalues only, namely, $|T(n,k)|$ with multiplicity $1$ and $-1$ with multiplicity $n! - 1$. Thus, for every $(n)\ne \zeta\vdash n$, formula (\ref{eq1}) evaluates to $-1$, that is, 
\begin{equation}\label{eq3}
\sum_{\substack{\gamma\vdash n \\ 2\leq n- c_1(\gamma)\leq n}} |C(S_n,\gamma)| \cdot \tilde{\chi}_\zeta(\gamma)=-1
\end{equation} 
for every $(n)\ne \zeta\vdash n$. So the maximum in formula (\ref{eq2}) also equals $-1$, which is attained by any $(n)\ne \zeta\vdash n$ and hence by $\zeta=(n-1,1)$ in particular. This means that the result in Theorem \ref{thm:result2} is true when $k=n$ with $n\ge 4$.

\smallskip
\textsf{Case 2.} $k = n-1$.
\smallskip

In this case, formula (\ref{eq2}) becomes \begin{equation*}
	\max_{\substack{\zeta\vdash n\\ \zeta\ne (n)}} \sum_{\substack{\gamma\vdash n\\ 2\leq n-c_1(\gamma)\leq n-1}} |C(S_n,\gamma)|\cdot \tilde{\chi}_{\zeta}(\gamma),
\end{equation*} which equals \begin{equation*}
	-1- \min_{\substack{\zeta\vdash n\\ \zeta\ne (n)}} \sum_{\substack{\gamma\vdash n \\c_1(\gamma)=0}} |C(S_n,\gamma)|\cdot \tilde{\chi}_{\zeta}(\gamma)
\end{equation*}
due to equation (\ref{eq3}). Thus, to prove that the second largest eigenvalue of $\mathrm{Cay}(S_n,T_{n-1})$ is given by the standard representation of $S_n$, it suffices to prove that the minimum of 
\begin{equation}\label{eq4}
	 \sum_{\substack{\gamma\vdash n\\ c_1(\gamma)=0}} |C(S_n,\gamma)|\cdot \tilde{\chi}_{\zeta}(\gamma)
\end{equation} 
among all $(n)\ne \zeta\vdash n$ is attained by $\zeta=(n-1,1)$. Note that, for $\zeta=(n-1,1)$ with $n\ge 4$ and any $\gamma \vdash n$ with $c_1(\gamma)=0$, we have $\tilde{\chi}_\zeta(\gamma)=\frac{c_1(\gamma)-1}{n-1}=-\frac{1}{n-1}$ by Table \ref{tab:tab1}. Hence the value of formula (\ref{eq4}) for $\zeta=(n-1,1)$ is 
\begin{equation}\label{eq5}
\frac{-1}{n-1} \cdot \sum_{\substack{\gamma\vdash n\\ c_1(\gamma)=0}} |C(S_n,\gamma)|.
\end{equation} 
In the following we will compare this with the value of (\ref{eq4}) for any $\zeta \vdash n$ with $\zeta \ne (n), (n-1,1)$.

By Lemma \ref{lem2.4}, if $\gamma$ is a partition of $n$ with $c_1(\gamma)=0$, then for every $\zeta\vdash n$ we have 
\begin{equation*} 
	\chi_\zeta(\gamma)\ge -\chi_\zeta(\iota)^{\frac{1}{2}+\varepsilon_n}
\end{equation*}  
and hence 
\begin{equation*}
	\tilde{\chi}_\zeta(\gamma)\ge -\chi_\zeta(\iota)^{-\frac{1}{2}+\varepsilon_n}.
\end{equation*}
Moreover, by Lemma \ref{lem2.5}, we have $\chi_\zeta(\iota)\ge n^{2.05}$ for any partition $\zeta$ of $n\ge 13$ whose Young diagram has at least three blocks outside the first row and at least three blocks outside the first column. Let $N_1$ be the smallest integer no less than $13$ such that $2.05 (-\frac{1}{2}+\varepsilon_n)\leq -1$ for all $n \ge N_1$, where $\varepsilon_n$ is as in Lemma \ref{lem2.4}. Then
\begin{equation*}
	\tilde{\chi}_\zeta(\gamma)\ge -\chi_\zeta(\iota)^{-\frac{1}{2}+\varepsilon_n}\ge -\frac{1}{n}> -\frac{1}{n-1}
\end{equation*}
for any $\zeta \vdash n$ satisfying the conditions in Lemma \ref{lem2.5}. So, when $n\ge N_1$, the values of (\ref{eq4}) corresponding to these $\zeta$'s are greater than (\ref{eq5}). Therefore, to determine the minimum of \eqref{eq4} it remains to consider the following five possibilities for $\zeta \ne (n-1,1)$: $(n-1, 1)'$, $(n- 2, 2)$, $(n- 2, 2)'$, $(n - 2, 1, 1)$, and $(n - 2, 1, 1)'$. Using Table \ref{tab:tab1} and Lemma \ref{lem2.1}, we obtain that the normalized character $\tilde{\chi}_\zeta(\gamma)$ of $\zeta = (n-1, 1)', (n- 2, 2), (n - 2, 1, 1), (n - 2, 1, 1)'$ on any $\gamma \vdash n$ with $c_1(\gamma)=0$ is equal to $\frac{-\mathrm{sgn}(\gamma)}{n-1}, \frac{2c_2(\gamma)}{n(n-3)}, \frac{2-2c_2(\gamma)}{(n-1)(n-2)}, \frac{\mathrm{sgn}(\gamma)(2-2c_2(\gamma))}{(n-1)(n-2)}$, respectively. Thus, if $n\ge 4$ and $\zeta = (n-1, 1)', (n- 2, 2), (n - 2, 1, 1), (n - 2, 1, 1)'$, then $\tilde{\chi}_\zeta(\gamma)\ge -\frac{1}{n-1}$ for any $\gamma\vdash n$ with $c_1(\gamma)=0$. Therefore, the value of (\ref{eq4}) corresponding to any one of these four $\zeta$'s are greater than (\ref{eq5}) when $n\ge 4$.  

Now we assume that $\zeta=(n-2,2)'$ and $\gamma$ is any partition of $n$ with $c_1(\gamma)=0$. We aim to prove that the value of (\ref{eq4}) for this $\zeta$ is also greater than (\ref{eq5}). In fact, by Table \ref{tab:tab1} and Lemma \ref{lem2.1}, we have
$$
\tilde{\chi}_\zeta(\gamma)=\mathrm{sgn}(\gamma)\cdot \tilde{\chi}_{(n-2,2)}(\gamma) =\mathrm{sgn}(\gamma)\cdot \frac{2c_2(\gamma)}{n(n-3)}.
$$ 
If $n\ge 4$ is odd, then we have $2c_2(\gamma)\leq n-3$ as $c_1(\gamma)=0$, which implies that $\tilde{\chi}_\zeta(\gamma)>-\frac{1}{n-1}$.  If $n\ge 4$ is even, then $2c_2(\gamma)=n$ or $2c_2(\gamma)\leq n-4$ due to $c_1(\gamma)=0$. If $2c_2(\gamma)\leq n-4$, then we also have $\tilde{\chi}_\zeta(\gamma)>-\frac{1}{n-1}$. If $2c_2(\gamma)=n$ (that is, $\gamma=(2^{\frac{n}{2}})$), then $\tilde{\chi}_\zeta(\gamma)=\mathrm{sgn}(\gamma)\cdot \frac{1}{n-3}$ is greater than $-\frac{1}{n-1}$ if and only if $\mathrm{sgn}(\gamma)=1$. Note that when $\gamma=(2^{\frac{n}{2}})$, $\mathrm{sgn}(\gamma)=(-1)^{\frac{n}{2}}=1$ if and only if $n\equiv 0 \pmod{4}$. Thus, for any $n\equiv 0 \pmod{4}$, we still have $\tilde{\chi}_\zeta(\gamma)>-\frac{1}{n-1}$ for any $\gamma\vdash n$ with $c_1(\gamma)=0$. Combining these, we obtain that the value of (\ref{eq4}) for $\zeta=(n-2,2)'$ is greater than (\ref{eq5}) whenever $n\ge 4$ and $n\not\equiv 2 \pmod{4}$. 

Now assume that $n\ge 4$ and $n\equiv 2 \pmod{4}$. Then only for $\gamma=(2^{\frac{n}{2}})$ is $\tilde{\chi}_\zeta(\gamma)=-\frac{1}{n-3}$ smaller than $-\frac{1}{n-1}$, and for any other partition $\gamma$ of $n$ with $c_1(\gamma)=0$ we still have $2c_2(\gamma)\leq n-4$ and hence $\tilde{\chi}_\zeta(\gamma)>-\frac{1}{n-1}$. Since $\mathrm{sgn}\left((4,2^{\frac{n-4}{2}})\right)=1$, by Table 1 we have
$$
\tilde{\chi}_{\zeta}((4,2^{\frac{n-4}{2}}))=\mathrm{sgn}\big((4,2^{\frac{n-4}{2}})\big)\cdot \tilde{\chi}_{(n-2,2)}\big((4,2^{\frac{n-4}{2}})\big)=\frac{n-4}{n(n-3)}.
$$ 
Thus,
\begin{eqnarray*}
& & \qquad \sum_{\substack{\gamma\vdash n\\ c_1(\gamma)=0}} |C(S_n,\gamma)| \cdot \tilde{\chi}_{(n-2,2)'}(\gamma) \\
		& = & \sum_{\substack{\gamma\vdash n,~c_1(\gamma)=0\\ \gamma \notin \{(2^{\frac{n}{2}}), (4,2^{\frac{n-4}{2}})\}}} |C(S_n,\gamma)| \cdot \tilde{\chi}_{(n-2,2)'}(\gamma) + |C(S_n,(2^{\frac{n}{2}}))|\cdot \tilde{\chi}_{(n-2,2)'}((2^\frac{n}{2})) \\ 
		& & \qquad \qquad \qquad \qquad +\ |C(S_n,(4,2^{\frac{n-4}{2}}))|\cdot \tilde{\chi}_{(n-2,2)'}((4,2^\frac{n-4}{2}))\\[2mm] 
		& > & \sum_{\substack{\gamma\vdash n,~ c_1(\gamma)=0\\ \gamma \notin \{(2^{\frac{n}{2}}), (4,2^{\frac{n-4}{2}})\}}} |C(S_n,\gamma)| \cdot \tilde{\chi}_{(n-1,1)}(\gamma)+|C(S_n,(2^{\frac{n}{2}}))|\cdot \tilde{\chi}_{(n-2,2)'}((2^\frac{n}{2}))\\[2mm]
		& & \qquad \qquad \qquad \qquad +\ |C(S_n,(4,2^{\frac{n-4}{2}}))|\cdot \tilde{\chi}_{(n-2,2)'}((4,2^\frac{n-4}{2}))\\[2mm]
		& = & \sum_{\substack{\gamma\vdash n,~ c_1(\gamma)=0\\ \gamma \notin \{(2^{\frac{n}{2}}), (4,2^{\frac{n-4}{2}})\}}} |C(S_n,\gamma)| \cdot \tilde{\chi}_{(n-1,1)}(\gamma)+|C(S_n,(2^{\frac{n}{2}}))|\cdot \frac{-1}{n-3}\\[2mm] 
		& & \qquad \qquad \qquad \qquad +\ |C(S_n,(4,2^{\frac{n-4}{2}}))|\cdot \frac{n-4}{n(n-3)}\\[2mm]
		& > & \sum_{\substack{\gamma\vdash n,~ c_1(\gamma)=0\\ \gamma \notin \{(2^{\frac{n}{2}}), (4,2^{\frac{n-4}{2}})\}}} |C(S_n,\gamma)| \cdot \tilde{\chi}_{(n-1,1)}(\gamma)\\[2mm]
		& > &\sum_{\substack{\gamma\vdash n,~ c_1(\gamma)=0}} |C(S_n,\gamma)| \cdot \tilde{\chi}_{(n-1,1)}(\gamma).
\end{eqnarray*}
In the second last step above, the inequality holds as $|C(S_n,(4,2^{\frac{n-4}{2}}))|=\frac{n(n-2)}{4}|C(S_n,(2^{\frac{n}{2}}))|$ and $|C(S_n,(2^{\frac{n}{2}}))|\cdot \frac{-1}{n-3}+|C(S_n,(4,2^{\frac{n-4}{2}}))|\cdot \frac{n-4}{n(n-3)}>0$. Thus, when $n\ge 4$ and $n\equiv 2 \pmod{4}$, the value of (\ref{eq4}) for $\zeta=(n-2,2)'$ is also greater than (\ref{eq5}). 

To sum up, we have proved that for any $n\ge N_1$ the minimum of \eqref{eq4} among all $(n)\ne \zeta\vdash n$ is attained by $\zeta=(n-1,1)$. In other words, the statement in Theorem \ref{thm:result2} is true for $k=n-1$ whenever $n\ge N_1$.

We claim that $N_1$ is no more than the integer $N_0$ in Lemma \ref{lem2.2}. In fact, by the definition of $N_1$, we know that $N_1$ is no more than $N_3$ in \cite[Lemma 2.7]{PP}, which is an integer satisfying $2.05 \left(-\frac{1}{2}+\frac{\log 2}{2 \log n}+\varepsilon_n\right) \leq-1$ for all $n\ge N_3$. Also, we see from the proof of \cite[Proposition 2.3]{PP} that $N_0$ is no less than $N_3$ in \cite[Lemma 2.7]{PP}. Hence $N_1 \le N_0$ as claimed.

\smallskip
\textsf{Case 3.} $2\leq k\leq n-2$.
\smallskip

This case is a direct consequence of Theorem \ref{thm:result1}. In fact, letting $I_k=\{2,3,\ldots,k\}$ for each $2\leq k\leq n-2$, we have $I_k \subseteq \{2,3,\ldots,n-2\}$ and $\mathrm{Cay}(S_n,T(n,k)) = \mathrm{Cay}(S_n, T(n,I_k))$ as $T(n,k) = T(n,I_k)$. So, by Theorem \ref{thm:result1} and its proof, there is a positive integer $N_2 := \max\{N_0, 7\}$ such that for every $n\ge N_2$ and any $2\leq k\leq n-2$, $\mathrm{Cay}(S_n,T(n,k))$ has the Aldous property, where $N_0$ is the integer in Lemma \ref{lem2.2}.
 
\vspace{0.2cm}
In summary, we have proved that for every $n\ge N := \max\{N_1,N_2\} = \max\{N_0,7\}$ and any $2\leq k\leq n$, $\mathrm{Cay}(S_n,T(n,k))$ has the Aldous property. 
\end{proof}

\begin{re}
\label{re:re2}
When dealing with the case $k=n-1$ in the proof above, the key was to prove the statement that the minimum of formula (\ref{eq4}) among all $(n)\ne \zeta\vdash n$ is attained by $\zeta=(n-1,1)$ when $n$ is sufficiently large. This statement is equivalent to the fact that the smallest eigenvalue of the derangement graph $\mathrm{Cay}(S_n,\mathcal{D}_n)$ is attained by the standard representation of $S_n$ when $n$ is sufficiently large, where $\mathcal{D}_n$ is the set of derangements on $[n]$. After completing the proof of Theorem \ref{thm:result2}, we realized that Renteln had proved a stronger result \cite[Theorem 7.1]{R}, which asserts that for $n \ge 4$ the smallest eigenvalue of $\mathrm{Cay}(S_n,\mathcal{D}_n)$ is equal to $-\frac{|\mathcal{D}_n|}{n-1}$, settling affirmatively a conjecture of Ku and Wong \cite[Conjecture 1]{KW}, and moreover for $n\ge 5$ this smallest eigenvalue is achieved uniquely by the standard representation of $S_n$. Renteln proved this result using a recurrence formula \cite[Theorem 6.5]{R}, while our proof above in the case $k=n-1$ took a different approach.
\end{re}

Finally, we prove Corollary \ref{cor1} with the help of Theorem \ref{thm:result1} and some results from \cite{DZ, KLW, R}.

\begin{proof}
In \cite{DZ}, Deng and Zhang proved that for $n\ge 4$ the second largest eigenvalue of $\mathcal{F}(n,0)$ is positive and is given by the irreducible representation of $S_n$ corresponding to the partition $(n-2,2)$ of $n$. Combining this and the fact \cite[Theorem 7.1]{R} that for $n \ge 5$ the smallest eigenvalue of $\mathcal{F}(n,0)$ is negative and is achieved by $\rho_{(n-1,1)}$ (see Remark \ref{re:re2}), we know that $\mathcal{F}(n,0)$ does not have the Aldous property when $n\ge 5$. 

In \cite{KLW},  Ku, Lau and Wong proved that for $n\ge 7$ the smallest eigenvalue of $\mathcal{F}(n,1)$ is achieved only by the irreducible representation of $S_n$ corresponding to the partition $(n-2,2)$. In \cite[Lemma 3.5]{KLW}, they also proved that for $n\ge 7$ the standard representation of $S_n$ yields the eigenvalue $0$ of $\mathcal{F}(n,1)$ while at least one of the partitions $(1^n),\ (2^2,1^{n-4}),\ (3,1^{n-3})$ produces a positive eigenvalue of $\mathcal{F}(n,1)$ other than $|T(n,\{n-1\})|$. This implies that the second largest eigenvalue of $\mathcal{F}(n,1)$ is not attained by the standard representation of $S_n$; that is, $\mathcal{F}(n,1)$ does not have the Aldous property when $n\ge 7$.

On the other hand, for $2 \leq k \leq n-2$, we have $\{n-k\} \subseteq \{2, 3, \ldots, n-2\}$. Thus, by Theorem \ref{thm:result1}, $\mathcal{F}(n,k) = \mathrm{Cay}(S_n, T(n,\{n-k\}))$ has the Aldous property whenever $n \ge N$, where $N$ is as in Theorem \ref{thm:result1}.
\end{proof}

\medskip
\noindent \textbf{Acknowledgements}
\medskip

We would like to thank the three anonymous referees for their helpful comments. The first author was supported by the Melbourne Research Scholarship provided by The University of Melbourne.

\end{document}